\pdfoutput=1
\RequirePackage{ifpdf}
\ifpdf 
\documentclass[pdftex]{sigma}
\else
\documentclass{sigma}
\fi

\numberwithin{equation}{section}

\newtheorem{Theorem}{Theorem}[section]

\newtheorem{Lemma}[Theorem]{Lemma}
\newtheorem{Proposition}[Theorem]{Proposition}
 { \theoremstyle{definition}

 }

\begin{document}

\newcommand{\arXivNumber}{1711.06009}

\renewcommand{\thefootnote}{}

\renewcommand{\PaperNumber}{033}

\FirstPageHeading

\ShortArticleName{The Duals of the 2-Modular Irreducible Modules of the Alternating Groups}

\ArticleName{The Duals of the 2-Modular Irreducible Modules\\ of the Alternating Groups\footnote{This paper is a~contribution to the Special Issue on the Representation Theory of the Symmetric Groups and Related Topics. The full collection is available at \href{https://www.emis.de/journals/SIGMA/symmetric-groups-2018.html}{https://www.emis.de/journals/SIGMA/symmetric-groups-2018.html}}}

\Author{John MURRAY}
\AuthorNameForHeading{J.~Murray}

\Address{Department of Mathematics \& Statistics, Maynooth University, Co. Kildare, Ireland}
\Email{\href{mailto:John.Murray@mu.ie}{John.Murray@mu.ie}}

\ArticleDates{Received January 04, 2018, in final form April 04, 2018; Published online April 17, 2018}

\Abstract{We determine the dual modules of all irreducible modules of alternating groups over fields of characteristic~2.}


\Classification{20C30}

\renewcommand{\thefootnote}{\arabic{footnote}}
\setcounter{footnote}{0}

\section{Introduction and statement of the result}

Let ${\mathcal S}_n$ be the symmetric group of degree $n\geq1$ and let $k$ be a field of characteristic $p>0$. In \cite[Theorem~11.5]{James} G.~James constructed all irreducible $k{\mathcal S}_n$-modules $D^\lambda$ where $\lambda$ ranges over the $p$-regular partitions of $n$. Here a partition is $p$-regular if each of its parts occurs with multiplicity less than~$p$.

As the alternating group ${\mathcal A}_n$ has index $2$ in ${\mathcal S}_n$, the restriction $D^\lambda{\downarrow_{{\mathcal A}_n}}$ is either irreducible or splits as a direct sum of two non-isomorphic irreducible $k{\mathcal A}_n$-modules. Moreover, every irreducible $k{\mathcal A}_n$-module is a direct summand of some $D^\lambda{\downarrow_{{\mathcal A}_n}}$.

Henceforth we will assume, unless stated otherwise, that $k$ is a field of characteristic $2$ which is a splitting field for the alternating group ${\mathcal A}_n$. For this, it suffices that $k$ contains the finite field ${\mathbb F}_4$. D.~Benson \cite{Benson} has classified all irreducible $k{\mathcal A}_n$-modules:

\begin{Proposition}\label{P:Benson}
Let $\lambda=({\lambda_1>\lambda_2>\dots>\lambda_{2s-1}>\lambda_{2s}\geq0})$ be a strict partition of $n$. Then $D^\lambda{\downarrow_{{\mathcal A}_n}}$ is reducible if and only if
\begin{itemize}\itemsep=0pt
 \item[$(i)$] $\lambda_{2j-1}-\lambda_{2j}=1$ or $2$, for $j=1,\dots,s$, and
 \item[$(ii)$] $\lambda_{2j-1}+\lambda_{2j}\not\equiv2$ $({\rm mod}~4)$, for $j=1,\dots,s$.
\end{itemize}
\end{Proposition}

In this note we determine the dual of each irreducible $k{\mathcal A}_n$-module. Now $D^\lambda{\downarrow_{{\mathcal A}_n}}$ is a self-dual $k{\mathcal A}_n$-module, as $D^\lambda$ is a self-dual $k{\mathcal S}_n$-module. So we only need to determine the dual of an irreducible $k{\mathcal A}_n$-module which is a direct summand of $D^\lambda{\downarrow_{{\mathcal A}_n}}$, when this module is reducible.

\begin{Theorem}\label{T:main}
Let $\lambda$ be a strict partition of $n$ such that $D^\lambda{\downarrow_{{\mathcal A}_n}}$ is reducible. Then the two irreducible direct summands of $D^\lambda{\!\downarrow_{{\mathcal A}_n}}$ are self-dual if ${\sum\limits_{j=1}^s \lambda_{2j}}$ is even and are dual to each other if ${\sum\limits_{j=1}^s \lambda_{2j}}$ is odd.
\end{Theorem}

For example $D^{(7,5,1)}{\downarrow_{{\mathcal A}_{13}}}\cong S\oplus S^*$, for a non self-dual irreducible $k{\mathcal A}_{13}$-module $S$, and $D^{(5,4,3,1)}{\downarrow_{{\mathcal A}_{13}}}$ decomposes similarly. On the other hand $D^{(7,6)}{\downarrow_{{\mathcal A}_{13}}}\cong S_1\oplus S_2$ where $S_1$ and $S_2$ are irreducible and self-dual.

In order to prove Theorem \ref{T:main}, we use the following elementary result, which requires the assumption that $k$ has characteristic $2$:

\begin{Lemma}\label{L:self-dual}
Let $G$ be a finite group and let $M$ be a semisimple $kG$-module which affords a~non-degenerate $G$-invariant symmetric bilinear form $B$. Suppose that $B(tm,m)\ne{ 0}$, for some involution $t\in G$ and some $m\in M$. Then $M$ has a self-dual irreducible direct summand.
\begin{proof}
We have $M=\bigoplus_{i=1}^n M_i$, for some $n\geq1$ and irreducible $kG$-modules $M_1,\dots,M_n$. Write $m=\sum m_i$, with $m_i\in M_i$, for all $i$. Then
\begin{gather*}
B(tm,m)=\sum_{1\leq i\leq n}B(tm_i,m_i)+\sum_{1\leq i<j\leq n}^n\big(B(tm_i,m_j)+B(tm_j,m_i)\big)\\
\hphantom{B(tm,m)}{} =\sum_{1\leq i\leq n}B(tm_i,m_i).
\end{gather*}
The last equality follows from the fact that $\operatorname{char}(k)=2$ and
\begin{gather*}
B(tm_i,m_j)=B\big(m_i,t^{-1}m_j\big)=B(m_i,tm_j)=B(tm_j,m_i).
\end{gather*}

Without loss of generality $B(tm_1,m_1)\ne0$. Then $B$ restricts to a non-zero $G$-invariant symmetric bilinear form $B_1$ on $M_1$. As $M_1$ is irreducible, $B_1$ is non-degenerate. So { $M_1$ is isomorphic to its $kG$-dual $M_1^*$}.
\end{proof}
\end{Lemma}

\section{Known results on the symmetric and alternating groups}

\subsection{The irreducible modules of the symmetric groups}\label{SS:irred}

We use the ideas and notation of \cite{James}. In particular for each partition $\lambda$ of $n$, James defines the Young diagram $[\lambda]$ of $\lambda$, and the notions of a $\lambda$-tableau and a $\lambda$-tabloid.

Fix a $\lambda$-tableau $x$. So $x$ is a filling of $[\lambda]$ with the symbols $\{1,\dots,n\}$. The corresponding $\lambda$-tabloid is $\{x\}:=\{\sigma(x)\,|\,\sigma\in R_x\}$, where $R_x$ is the row stabilizer of $x$. We regard $\{x\}$ as an ordered set partition of $\{1,\dots,n\}$. The ${\mathbb Z}$-span of the $\lambda$-tabloids forms the ${\mathbb Z}{\mathcal S}_n$-lattice $M^\lambda$, and the set of $\lambda$-tabloids is an ${\mathcal S}_n$-invariant ${\mathbb Z}$-basis of $M^\lambda$.

Recall from \cite[Section~4]{James} that corresponding to each tableau $x$ there is a polytabloid $e_x:=\sum\operatorname{sgn}(\sigma)\{\sigma x\}$ in $M^\lambda$. Here $\sigma$ ranges over the permutations in the column stabilizer $C_x$ of the tableau $x$. The Specht lattice $S^\lambda$ is defined to be the ${\mathbb Z}$-span of all $\lambda$-polytabloids. In particular~$S^\lambda$ is a ${\mathbb Z}{\mathcal S}_n$-sublattice of $M^\lambda$; it has as ${\mathbb Z}$-basis the polytabloids corresponding to the standard $\lambda$-tableaux (i.e., the numbers increase from left-to-right along rows, and from top-to-bottom along columns).

Now James defines $\langle\,,\,\rangle$ to be the symmetric bilinear form on $M^\lambda$ which makes the tabloids into an orthonormal basis. As the tabloids are permuted by the action of ${\mathcal S}_n$, it is clear that $\langle\,,\,\rangle$ is ${\mathcal S}_n$-invariant.

Suppose now that $\lambda$ is a strict partition and consider the unique permutation $\tau\in R_x$ which reverses the order of the symbols in each row of the tableau $x$. In \cite[Lemma~10.4]{James} James shows that $\langle \tau e_x,e_x\rangle=1$, as $\{x\}$ is the only tabloid common to $e_x$ and $e_{\tau x}$ (in fact James proves that $\langle \tau e_x,e_x\rangle$ is coprime to $p$, if $\lambda$ is $p$-regular, for some prime $p$). Set $J^\lambda:=\{x\in S^\lambda\,|\, \langle x,y\rangle\in2{\mathbb Z}$, for all \smash{$y\in S^\lambda\}$}. Then $2S^\lambda\subseteq J^\lambda$ and it follows from \cite[Theorem~4.9]{James} that $D^\lambda:=(S^\lambda/J^\lambda)\otimes_{\mathbb F_2}k$ is an absolutely irreducible $k{\mathcal S}_n$-module, for any field $k$ of characteristic~$2$.

\subsection{The real 2-regular conjugacy classes of the alternating groups}

A conjugacy class of a finite group $G$ is said to be $2$-regular if its elements have odd order. R.~Brauer proved that the number of irreducible $kG$-modules equals the number of $2$-regular conjugacy classes of~$G$~\cite{BrauerNesbitt}. Now Brauer's permutation lemma holds for arbitrary fields \cite[footnote~19]{Brauer}. So it is clear that the number of self-dual irreducible $kG$-modules equals the number of real $2$-regular conjugacy classes of~$G$.

We review some well known facts about the $2$-regular conjugacy classes of the alternating group. See for example \cite[Section~2.5]{JamesKerber}.

Corresponding to each partition $\mu$ of $n$ there is a conjugacy class $C_\mu$ of ${\mathcal S}_n$; its elements consist of all permutations of $n$ whose orbits on $\{1,\dots,n\}$ have sizes $\{\mu_1,\dots,\mu_\ell\}$ (as multiset). So~$C_\mu$ is $2$-regular if and only if each $\mu_i$ is odd.

Let $\mu$ be a partition of $n$ into odd parts. Then $C_\mu\subseteq{\mathcal A}_n$. If $\mu$ has repeated parts then $C_\mu$ is a conjugacy class of ${\mathcal A}_n$. As $C_\mu$ is closed under taking inverses, $C_\mu$ is a real conjugacy class of~${\mathcal A}_n$.

Now assume that $\mu$ has distinct parts. Then $C_\mu$ is a union of two conjugacy classes $C_\mu^{\pm}$ of~${\mathcal A}_n$. Set $m:=\frac{n-\ell(\mu)}{2}$ and let $z\in C_\mu$. Then $z$ is inverted by an involution $t\in{\mathcal S}_n$ of cycle type $(2^m,1^{n-2m})$. Since $\operatorname{C}_{{\mathcal S}_n}(z)\cong\prod {\mathbb Z}/{\mu_j}{\mathbb Z}$ is odd, $t$ generates a Sylow $2$-subgroup of the extended centralizer $\operatorname{C}_{{\mathcal S}_n}^*(z)$ of $z$ in ${\mathcal S}_n$. It follows that $z$ is conjugate to $z^{-1}$ in ${\mathcal A}_n$ if and only if $t\in{\mathcal A}_n$. This shows that $C_\mu^{\pm}$ are real classes of ${\mathcal A}_n$ if and only if $\frac{n-\ell(\mu)}{2}$ is even. This and the discussion above shows:

\begin{Lemma}\label{L:non-self-dual}
The number of self-dual irreducible $k{\mathcal A}_n$-modules equals the number of non-strict odd partitions of $n$ plus twice the number of strict odd partitions $\mu$ of $n$ for which $\frac{n-\ell(\mu)}{2}$ is even.
\end{Lemma}

\section{Bressoud's bijection}

We need a special case of a partition identity of I.~Schur \cite{Schur}. This was already used by Benson in his proof of Proposition~\ref{P:Benson}:

\begin{Proposition}[Schur, 1926]\label{P:Schur}
The number of strict partitions of $n$ into odd parts equals the number of strict partitions of $n$ into parts congruent to $0$, $\pm1$ $({\rm mod}~4)$ where consecutive parts differ by at least $4$ and consecutive even parts differ by at least $8$.
\end{Proposition}

D.~Bressoud \cite{Bressoud} has constructed a bijection between the relevant sets of partitions. We describe a simplified version of this bijection.

Let $\mu=(\mu_1>\mu_2>\dots>\mu_\ell)$ be a strict partition of $n$ whose parts are all odd. We sub\-divide~$\mu$ into `blocks' of at most two parts, working recursively from largest to smallest parts. Let $j\geq1$ and suppose that $\mu_1,\mu_2,\dots,\mu_{j-1}$ have already been assigned to blocks. We form the block $\{\mu_j,\mu_{j+1}\}$ if $\mu_j=\mu_{j+1}+2$, and the block $\{\mu_j\}$ otherwise (if $\mu_j\geq\mu_{j+1}+4$). Let $s$ be the number of resulting blocks of $\mu$.

Next we form the sequence of positive integers $\sigma=(\sigma_1,\sigma_2,\dots,\sigma_s)$, where $\sigma_j$ is the sum of the parts in the $j$-th block of $\mu$. Then the $\sigma_j$ are distinct, as the odd parts form a decreasing sequence, with minimal difference $4$, and the even parts form a decreasing sequence, with minimal difference $8$. Moreover, each even $\sigma_j$ is the sum of a pair of consecutive odd integers. So $\sigma_j\not\equiv2$ $({\rm mod}~4)$, for all $j>0$.

We get a composition $\zeta$ of $n+2s(s-1)$ by defining
\begin{gather*}
\zeta_1=\sigma_1,\ \zeta_2=\sigma_2+4,\ \dots, \ \zeta_s=\sigma_s+4(s-1).
\end{gather*}
The even $\zeta_j$ form a decreasing sequence, with minimal difference $4$, and the odd $\zeta_j$ form a weakly decreasing sequence ($\zeta_j=\zeta_{j+1}$ if and only if $\zeta_j$, $\zeta_{j+1}$ represent two singleton blocks $\{2k-1\}$ and $\{2k-5\}$ of $\mu$, for some $k\geq0$).

Choose a permutation $\tau$ such that $\zeta_{\tau1}\geq\zeta_{\tau2}\geq\dots\geq\zeta_{\tau s}$. { Then we get a strict partition $\gamma$ of $n$ by defining
\begin{gather*}
\gamma_1=\zeta_{\tau1},\ \gamma_2=\zeta_{\tau2}-4, \ \dots, \ \gamma_s=\zeta_{\tau s}-4(s-1).
\end{gather*}
By construction, the minimal difference between the parts of $\gamma$ is $4$ and the minimal difference between the even parts of $\gamma$ is $8$. Moreover, $\gamma_j\equiv\zeta_{\tau j}$ $({\rm mod}~4)$. So $\gamma_j\not\equiv2$ $({\rm mod}~4)$. Then $\mu\rightarrow\gamma$ is Bressoud's bijection.

Finally form a strict partition $\lambda$ of $n$ which has $2s-1$ or $2s$ parts, by defining
\begin{gather*}
(\lambda_{2j-1},\lambda_{2j})=
 \begin{cases}
\displaystyle \left(\frac{\gamma_j}{2}+1,\frac{\gamma_j}{2}-1\right),& \mbox{if $\gamma_j$ is even or}\vspace{1mm}\\
\displaystyle \left(\frac{\gamma_j+1}{2},\frac{\gamma_j-1}{2}\right),& \mbox{if $\gamma_j$ is odd}.
 \end{cases}
\end{gather*}
Then $\lambda$ satisfies the constraints (i) and (ii) of Proposition~\ref{P:Benson}. Conversely, it is easy to see that if $\lambda$ satisfies these constraints, then $\lambda$ is the image of some strict odd partition $\mu$ of $n$ under the above sequence of operations.

\begin{Lemma}\label{L:n-l=sum}
Let $\mu$ be a strict-odd partition of $n$ and let $\lambda$ be the strict partition of $n$ constructed from $\mu$ as above. Then $\frac{n-\ell(\mu)}{2}=\sum\lambda_{2j}$.
\begin{proof}
Each pair of consecutive parts $\lambda_{2j-1}$, $\lambda_{2j}$ of $\lambda$ corresponds to a block ${\mathcal B}$ of~$\mu$. Moreover by our description of Bressoud's bijection, there are integers $q_1,\dots,q_s$, with $\sum q_j=0$ such that
\begin{gather*}
(\lambda_{2j-1}+2q_j,\lambda_{2j}+2q_j)=
\begin{cases}
\displaystyle \left(\frac{\mu_i+1}{2},\frac{\mu_i-1}{2}\right),& \mbox{if ${\mathcal B}=\{\mu_i\}$},\\
(\mu_i,\mu_{i+1}),& \mbox{if ${\mathcal B}=\{\mu_i,\mu_{i+1}\}$}.
\end{cases}
\end{gather*}
In case ${\mathcal B}=\{\mu_i,\mu_{i+1}\}$, we have $\mu_i=\mu_{i+1}+2$ and thus $\frac{\mu_i-1}{2}+\frac{\mu_{i+1}-1}{2}=\lambda_{2j}+2q_j$. We conclude that
\begin{gather*}
\frac{n-\ell(\mu)}{2}=\sum_{i=1}^{\ell(\mu)}\frac{\mu_i-1}{2}=\sum_{j=1}^s(\lambda_{2j}+2q_j)=\sum_{j=1}^s\lambda_{2j}.\tag*{\qed}
\end{gather*}\renewcommand{\qed}{}
\end{proof}
\end{Lemma}

\section{Proof of Theorem \ref{T:main}}

Let $\operatorname{D}(n)$ be the set of strict partitions of $n$ and let $\operatorname{S}(n)$ be the set of strict partitions of $n$ which satisfy conditions (i) and (ii) in Proposition \ref{P:Benson}. So there are $2|\operatorname{S}(n)|+|\operatorname{D}(n)\backslash\operatorname{S}(n)|$ irreducible $k{\mathcal A}_n$-modules.

Next set $\operatorname{S}(n)^+:=\{\lambda\in\operatorname{S}(n)\,|\,\sum\mbox{$\lambda_{2j}$ is even}\}$. Then it follows from Lemmas~\ref{L:non-self-dual} and~\ref{L:n-l=sum} that the number of self-dual irreducible $k{\mathcal A}_n$-modules equals $2|\operatorname{S}(n)^+|+|\operatorname{D}(n)\backslash\operatorname{S}(n)|$. Now $D^\lambda{\downarrow_{{\mathcal A}_n}}$ is an irreducible self-dual $k{\mathcal A}_n$-module, for $\lambda\in\operatorname{D}(n)\backslash\operatorname{S}(n)$. So we can prove Theorem~\ref{T:main} by showing that the irreducible direct summands of $D^\lambda{\downarrow_{{\mathcal A}_n}}$ are self-dual for all $\lambda\in\operatorname{S}(n)^+$.

Suppose then that $\lambda\in\operatorname{S}(n)^+$. Let $\tau\in{\mathcal S}_n$ be the permutation which reverses each row of a~$\lambda$-tableau, as discussed in Section~\ref{SS:irred}. We claim that $\tau\in{\mathcal A}_n$. For $\tau$ is a product of $\sum\limits_{i=1}^{2s}\big\lfloor\frac{\lambda_j}{2}\big\rfloor$ commuting transpositions. Now $\big\lfloor\frac{\lambda_{2j-1}}{2}\big\rfloor+\big\lfloor\frac{\lambda_{2j}}{2}\big\rfloor=\lambda_{2j}$, as $\lambda_{2j-1}-\lambda_{2j}=1$, or $\lambda_{2j-1}-\lambda_{2j}=2$ and both $\lambda_{2j-1}$ and $\lambda_{2j}$ are odd. So $\sum\limits_{i=1}^{2s}\big\lfloor\frac{\lambda_i}{2}\big\rfloor=\sum\limits_{j=1}^s\lambda_{2j}$ is even. This proves the claim.

Since $D^\lambda$ is irreducible and the form $\langle\,,\,\rangle$ is non-zero, $\langle\,,\,\rangle$ is non-degenerate on $D^\lambda$. Write $D^\lambda{\downarrow_{{\mathcal A}_n}}=S_1\oplus S_2$, where $S_1$ and $S_2$ are non-isomorphic irreducible modules. As $\tau\in{\mathcal A}_n$, it follows from Lemma \ref{L:self-dual} that we may assume that $S_1$ is self-dual. Now $S_2^*\not\cong S_1^*\cong S_1$ and $S_2^*$ is isomorphic to a direct summand of $D^\lambda{\downarrow_{{\mathcal A}_n}}$. So $S_2$ is also self-dual. This completes the proof of the theorem.

\section[Irreducible modules of alternating groups over fields of odd characteristic]{Irreducible modules of alternating groups\\ over fields of odd characteristic}

We now comment briefly on what happens when $k$ is a splitting field for ${\mathcal A}_n$ which has odd characteristic $p$. Let $\operatorname{sgn}$ be the sign representation of $k{\mathcal S}_n$. So $\operatorname{sgn}$ is $1$-dimensional but non-trivial. G.~Mullineux defined a bijection $\lambda\rightarrow\lambda^M$ on the $p$-regular partitions of $n$ and conjectured that $D^\lambda\otimes\operatorname{sgn}=D^{\lambda^M}$ for all $p$-regular partitions $\lambda$ of $n$. This was only proved in the 1990's by Kleshchev and Ford--Kleshchev. See~\cite{FordKleschev} for details.

Now $D^\lambda{\!\downarrow_{{\mathcal A}_n}}\cong D^{\lambda^M}{\!\downarrow_{{\mathcal A}_n}}$, and $D^\lambda{\!\downarrow_{{\mathcal A}_n}}$ is irreducible if and only if $\lambda\ne\lambda^M$ See \cite{Bessenrodt} for details. Moreover $D^\lambda$ and $D^{\lambda^M}$ are duals of each other, by \cite[Theorem~6.6]{James}. So $D^\lambda{\!\downarrow_{{\mathcal A}_n}}$ is self-dual, if $\lambda\ne\lambda^M$. However when $\lambda=\lambda^M$, we do not know how to determine when the two irreducible direct summands of $D^\lambda{\!\downarrow_{{\mathcal A}_n}}$ are self-dual.

\subsection*{Acknowledgement}

D.~Benson told me that it was an open problem to determine the self-dual irreducible $k{\mathcal A}_n$-modules. G.E.~Andrews directed me to Bressoud's paper~\cite{Bressoud}. We also thank the anonymous referees for their comments, which helped to improve the clarity of this paper.

\pdfbookmark[1]{References}{ref}
\LastPageEnding

\end{document}